\documentclass[12pt,centertags,oneside]{amsart}

\usepackage{dsfont}
\usepackage{amsmath,amstext,amsthm,amscd,typearea,hyperref,stmaryrd}
\usepackage{amssymb}
\usepackage{a4wide}
\usepackage[mathscr]{eucal}
\usepackage{mathrsfs}
\usepackage{typearea}
\usepackage{charter}
\usepackage{pdfsync}
\usepackage[a4paper,width=16.2cm,top=3cm,bottom=3cm]{geometry}

\numberwithin{equation}{section}

\newtheorem{theorem}{Theorem}[section]

\newtheorem{definition}[theorem]{Definition}
\newtheorem{proposition}[theorem]{Proposition}

\newtheorem{lemma}[theorem]{Lemma}

\newcommand\pro{\noindent{\em{Proof}.\ }}

\newcommand{\id}{{\rm id}}

\renewcommand{\Re}{{\rm Re}}

\newcommand{\C}{\mathbb{C}}
\newcommand{\D}{\mathbb{D}}

\newcommand{\N}{\mathbb{N}}

\newcommand{\R}{\mathbb{R}}

\newcommand{\B}{\mathbb{B}}

\title{Zalcman's  renormalization lemma,
Pinchuk's rescaling method, and Catlin's  estimates revisited}

\title[Zalcman's  lemma, Pinchuk's rescaling method, and Catlin's  estimates revisited]{Zalcman's  renormalization lemma, Pinchuk's rescaling method, and Catlin's  estimates revisited}

\author{Fran\c cois Berteloot}
\address{Universit\'e Toulouse 3, 
Institut Math\'ematique de Toulouse, 
118 route de Narbonne,
F-31062 Toulouse Cedex 9, France. }
\email{francois.berteloot@math.univ-toulouse.fr}

\date{}

\begin{document}

\maketitle

\begin{abstract}  We present  a  renormalization lemma  for certain maps defined on the unit disc of $\C$ and taking values in some metric space. We show that  
 the classical renormalization lemmas of  Zalcman and Miniowitz can be deduced from our lemma. We also use it to establish a general normality statement  for the Pinchuk's scaling method in $\C^2$ and, incidentally, reprove the Catlin's estimates for the Kobayashi metric in finite type domains.

\end{abstract}

\section{Introduction and main result}

"Nihil est in infinito quod non prius fuerit  in finito". This aphorism, coined by Andr\'e Bloch in his 1926 article \cite{Bloch}, was supposed to describe a general principle in complex
analysis, nowadays known as Bloch's principle.
In a 1973 address \cite{Rob}, Abraham Robinson  formulated a mathematically precise version of such a principle, and outlined a possible proof based on non standard analysis. 
Two years later, Lawrence Zalcman gave a standard, and very elementary, proof of the following lemma \cite{Zal}.

\begin{lemma}\label{Lemzal} Let $(f_n)_n$ be a sequence of holomorphic maps from the unit disc $\D$ of $\C$ to the Riemann sphere  which is not normal at $0$.
Then there exists a sequence of affine contractions $(r_n)_n$, converging to $0$ and such that the renormalized sequence $(f_n\circ r_n)_n$ is converging to a non constant  entire map $\varphi$, after taking a subsequence. Moreover, the spherical derivative of $\varphi$ is uniformly bounded:$\vert \varphi '\vert _{\sigma} \le \vert \varphi ' (0)\vert _{\sigma}=1$.
\end{lemma}

 According to this result, if a property $\mathcal P$ of meromorphic functions is stable both  by affine change of coordinates at the source and local uniform convergence, then
 very few meromorphic functions on $\D$ satisfy $\mathcal P$ if there are no non constant entire meromorphic functions which satisfy it:
\begin{eqnarray*}
&\;& \{f :\; \textrm{meromorphic on}\; \C\;\textrm{and  satisfying}\; {\mathcal P}\} \subset {\C\cup\{\infty\}}\;\;\;\;
 \Rightarrow\\ 
&\;& \{f :\; \textrm{meromorphic on}\; \D\;\textrm{and satisfying}\; {\mathcal P}\} \;\textrm{is normal}.
 \end{eqnarray*}
 This implication actually establishes the Robinson version of Bloch's principle.  For more material on 
 Bloch's principle and  a detailed description of its applications, we refer to Walter Bergweiler's review article \cite{B}.\\
 
 Our aim is to prove a quite general  renormalization lemma, working for maps defined on the unit disc of $\C$, with values in a metric space, and satisfying some Schwarz type property. As we shall see, this lemma covers Zalcman's result and its extension to quasi-regular maps due to Miniowitz. We shall also show that it
leads to a general normality statement for the Pinchuk's rescaling method. This new approach considerably
simplifies  our previous one developped in \cite{Be}, it relies neither on the existence of good peak functions nor on estimates of invariant metrics, but only on the basic geometric properties of finite type hypersurfaces, via an asymptotic Bloch principle.

We shall work with metric spaces  $(X,d)$  endowed with specific families $(J_{\eta})_{\eta\in X}$ of positive real valued functions.
We have  two main examples  in mind. On one hand  $J_\eta(\eta')=d(\eta,\eta')$ where $X$ is a compact manifold with a metric $d$, which corresponds to the Zalcman or Miniowitz cases, and $J_{\eta}(\eta')=\Vert A_\eta(\eta')\Vert$
where $X\subset \C^2$ and $(A_{\eta})_{\eta\in X}$ is a family of holomorphic automorphisms of $\C^2$ on the other hand, which corresponds to the Pinchuk's rescaling method.

Before stating our result, we must define the Schwarz property on which it is based. We denote by $\D_r$ the disc of radius $r$ centred at the origin in $\C$.

\begin{definition}
Let $(J_{\eta})_{\eta\in X}$ be a family of of positive functions on a metric space $(X,d)$ such that $J_\eta(\eta)=0$ for every $\eta\in X$.
Let $0<\alpha^-<\alpha^+ < 1$ and $c>0$ be some constants, and $s:[0,\alpha^+[ \to \R^+$ be a function which is  vanishing and continuous at $0$. We shall say that a map $f: \D_\rho\to X$ satisfies the \emph{Schwarz-type property} $S(J,\alpha^{\pm},c, s)$ on $\D_\rho$  if the following estimates occur for any $t_0,\epsilon_0 \in \C$ such that $\vert t_0\vert +\vert \epsilon_0\vert <\rho$:
$$\sup_{t\in \D_{\alpha^+}} J_{f(t_0)}(f(t_0 +t\epsilon_0)) \le c \;\Rightarrow\; J_{f(t_0)}(f(t_0 +t\epsilon_0)) \le s(\vert t\vert),\;\;\forall t\in \D_{\alpha^-}.$$
\end{definition}

It might be useful to note that if $f: \D\to X$ satisfies the Schwarz-type property $S(J,\alpha^{\pm},c, s)$ on $\D$ and $0<\rho$, then $t\mapsto f(\rho t)$ satisfies it on $\D_{\frac{1}{\rho}}$.\\

 We may now state our generalized version of the Zalcman Lemma.

\begin{lemma}

\label{NewZalc} Let $(X,d)$ be a metric space and $(J_{\eta})_{\eta\in X}$ be a family of positive real valued functions on $X$
such that: 
 \begin{itemize}
\item[i)] $J_\eta(\eta)=0$ for every $\eta \in X$,
\item[ii)]  $\lim_{\tau \to 0} \sup_{\eta \in K} \left(\sup_{d(\eta,\eta')\le \tau} J_\eta(\eta')\right)=0$ for every compact subset $K$ of $X$.
\end{itemize}
Let $f_n : {\overline \D}\to X$ be a sequence of continuous maps such that:
\begin{itemize}
\item[1-] $f_n$ satisfies the Schwarz-type property $S(J,\alpha^{\pm},c, s)$ on $\D$  for every $n\in \N$,
\item[2-] there exist a constant $k\in ]0,1]$ and sequences $(t'_n)_n$, $(\epsilon'_n)_n$ in $\C$ such that $\lim_n \epsilon'_n=0$, $2\vert t'_n\vert +\vert\epsilon'_n\vert <1$ and $J_{f_n(t'_n)}\left(f_n(t'_n+\epsilon'_n)\right)\ge k$ for every $n\in \N$. 
\end{itemize} 
Then there exists  a sequence of affine contractions $\left(r_n(t):=t_n+\epsilon_n t)\right)_n$ and  a sequence of positive real numbers $(R_n)_n$ such that $\lim_n\epsilon_n=0$, 
 $\lim_nR_n=+\infty$, $\vert t_n\vert +\vert\epsilon_n\vert <1$  for every $n\in \N$ and 
\begin{itemize}
\item[1'-] $g_n:=f_n\circ r_n$ is defined on $\D_{R_n+\alpha^+}$,
\item[2'-]  $J_{g_n(0)}\left(g_n(1)\right)\ge k$ and $J_{g_n(t)}\left(g_n(t+u)\right)\le s(\vert u\vert)$ for every $n\in \N$, every $t\in  D_{R_n}$ and every $u\in \D_{\alpha^-}$.
\end{itemize} 
If moreover $\lim_n t'_n =0$, then the contractions $r_n$ can be chosen so that  $\lim_n r_n(0) =0$.
\end{lemma}

\pro We may assume that $k<c$. Let us set $\rho_n:=\vert \epsilon'_n\vert$, then $\rho_n>0$ and $\lim_n\rho_n=0$.
For every $n\in \N$, we define
$${\mathcal D_n}:=\{(t,\epsilon)\in \C\times \C\;\colon\; \vert t\vert +\vert \epsilon\vert <1\;\textrm{and}\;  J_{f_n(t)}\left(f_n(t+\epsilon)\right)\ge k\}$$
$$\sigma_n:=\inf_{(t,\epsilon)\in {\mathcal D}_n} \frac{\vert\epsilon\vert}{1-\vert t\vert}.$$
As $\sigma_n \ge \inf_{(t,\epsilon)\in {\mathcal D}_n} \vert\epsilon\vert$ and $f_n({\overline  \D})=:K$ is compact in $(X,d)$, we must have $\sigma_n>0$
otherwise we would find $(f_n(t),f_n(t+\epsilon))=:(\eta,\eta')\in K\times X$ with $d(\eta,\eta')$  arbitrarily small and $J_\eta(\eta')\ge k$ which, by our assumption on $J$, is impossible. Moreover, as  $(t'_n,\epsilon'_n)\in {\mathcal D}_n$ and $2\vert t'_n\vert +\vert \epsilon'_n\vert <1$, we have
\begin{eqnarray}\label{sigma}
0<\sigma_n <\frac{\rho_n}{1-\frac{1}{2}}=2\rho_n.
\end{eqnarray}
Let us now fix $\alpha$ such that $0<\alpha^-<\alpha^+<\alpha<1$. After removing the first terms of the sequence $(f_n)_n$, we may assume that
\begin{eqnarray}\label{sigma/alpha}
\sqrt{\sigma_n}+\alpha <\frac{\alpha}{\sqrt{\sigma_n}}\;\;\textrm{and}\;\; \alpha (1-\sqrt{\sigma_n}) >\alpha^+,\; \forall n\in \N.
\end{eqnarray}
Let us now pick $(t_n,\epsilon_n)\in {\mathcal D}_n$ such that
\begin{eqnarray}\label{tepsilon}
\sigma_n \le \frac{\vert \epsilon_n\vert}{1-\vert t_n\vert} < \frac{\sigma_n}{\alpha},
\end{eqnarray}
and set $r_n(t):=t_n+t\epsilon_n$, $g_n:=f_n\circ r_n$ and $R_n:=\frac{\alpha}{\sqrt{\sigma_n}}$, for every $n\in \N$.
From (\ref{tepsilon}) and  (\ref{sigma/alpha}) one gets
$\frac{1-\vert t_n\vert}{\vert \epsilon_n\vert}> \frac{\alpha}{\sigma_n} > \frac{\alpha}{\sqrt{\sigma_n}}+1=R_n+1$, which shows that
the map $g_n$ is defined on $\D_{R_n+ 1}$.\\

By construction $J_{g_n(0)}\left(g_n(1)\right)\ge k$. Let us show that $J_{g_n(t)}\left(g_n(t+u)\right)\le s(\vert u\vert)$ for every $(t,u)\in  D_{R_n}\times \D_{\alpha^-}$. We proceed by contradiction and assume that there exists $\tau_0$ and $u'_0$ in $\C$ such that $\vert \tau_0 \vert < R_n$, $\vert u'_0\vert < \alpha^-$ and $J_{g_n(\tau_0)} (g_n(\tau_0 + u'_0)) > s(\vert u'_0\vert)$. Setting $\tilde{t}_n:=t_n+\epsilon_n \tau_0$, this can be rewritten as
\begin{eqnarray}\label{contra}
J_{f_n(\tilde{t}_n)} (f_n(\tilde{t}_n +\epsilon_n u'_0)) > s(\vert u'_0\vert).
\end{eqnarray}
Observe that (\ref{sigma/alpha}) yields $1+\vert \tau_0\vert <1+R_n= 1+\frac{\alpha}{\sqrt{\sigma_n}} < \frac{\alpha}{\sigma_n}$. It thus follows that 
$$1-\vert \tilde{t}_n\vert -\vert \epsilon_n\vert > 1-\vert t_n\vert 
-\vert \epsilon_n\vert(1+ \vert \tau_0\vert)> 1-\vert t_n\vert -\frac{\alpha\vert \epsilon_n\vert}{\sigma_n} =\vert \epsilon_n\vert\left( \frac{1-\vert t_n\vert}{\vert \epsilon_n\vert}- \frac{\alpha}{\sigma_n}\right) > 0,$$
where the last inequality is given by (\ref{tepsilon}). Now, as $f_n$ satisfies the Schwarz-type property $S(J,\alpha^{\pm},c, s)$ and
$1 >\vert \tilde{t}_n\vert +\vert \epsilon_n\vert$, we deduce from (\ref{contra}) that there exists $u_0$ such that 
\begin{eqnarray}\label{contra2}
u_0\in \D_{\alpha^+}\;\;\textrm{and}\;\;J_{f_n(\tilde{t}_n)} (f_n(\tilde{t}_n +\epsilon_n u_0)) > c>k.
\end{eqnarray}
This means that $(\tilde{t}_n,\epsilon_n u_0)\in {\mathcal D}_n$. Therefore $\frac{\vert\epsilon_n u_0\vert}{1-\vert \tilde{t}_n\vert} \ge \sigma_n$ and then
\begin{eqnarray*}
\vert u_0\vert \ge \frac{1-\vert \tilde{t}_n\vert}{\vert\epsilon_n\vert} \sigma_n  \ge
  \frac{1-\vert \tilde{t}_n\vert}{1-\vert {t}_n\vert}\alpha  &\ge&
  \frac{1-\vert {t}_n\vert - \vert \epsilon_n\vert \vert \tau_0\vert}{1-\vert {t}_n\vert}\alpha \ge 
\left(1- \frac{\vert \epsilon_n\vert \vert \tau_0\vert}{1-\vert {t}_n\vert}\right)\alpha\\
 &\ge& 
\left(1- \frac{\sigma_n}{\alpha}\vert \tau_0\vert\right)\alpha \ge 
(1-\sqrt{\sigma_n})\alpha
\end{eqnarray*}
where the second and fifth  inequalities come from (\ref{tepsilon}), and the last one comes from $\tau_0 < R_n= \frac{\alpha}{\sqrt{\sigma_n}}$.
By  (\ref{sigma/alpha}) this yields $\vert u_0\vert > \alpha^+$, 
which contradicts (\ref{contra2}).

Finally, when $\lim_n t'_n=0$, we may take a sequence of positive numbers $(\rho_n)_n$ which is converging sufficently slowly  to $0$ so that $\lim_n \frac{\epsilon'_n}{\rho_n}=0$ and $2\vert\frac{t'_n}{\rho_n}\vert + \vert\frac{\epsilon'_n}{\rho_n}\vert<1$
for every $n$. Applying the above result to the sequence $(\tilde{f}_n)_n$ where $\tilde{f}_n(t):=f_n(\rho_n t)$ yields contractions of the form  $\rho_n(t_n+\epsilon_n t)$ which renormalize 
$(f_n)_n$ and satisfy $\lim_n \rho_n t_n =0$.
\qed\\

To end this section let us stress that  the above proof clearly remains valid if the complex plane $\C$
is replaced by $\R^m$ or, \emph{mutatis mutandis}, by any homogeneous group like, for instance, the Heisenberg group.

\section{Zalcman and Miniowitz lemmas}

For suitable families of functions $(J_\eta)_{\eta\in X}$, the conclusion $2'$ in Lemma \ref{NewZalc} both guarantees that the renormalized sequence $(g_n)_n$ is equicontinuous on compact subsets of $\C$ and that 
its limits $g$ satisfy $J_{g(0)}(g(1)) \ge k$ and are therefore non constants. Since holomorphic, or quasi-regular maps, satisfy a Schwarz lemma, this facts explain why  Zalcman's and Miniowitz's results can be easily derived from our general renormalization lemma.\\

Let us now give some details. We start with the case of holomorphic maps from the unit disc to some compact complex manifold $X$ of dimension $m\ge 1$, a setting which is slightly more general than that of meromorphic functions considered in Lemma \ref{Lemzal}. Let $d$ be a distance induced by an hermitian metric on $X$ and let  $B(\eta,r)$ denote the open ball centred at $\eta$ and of radius $r>0$ relative to this distance. Let $\Vert\;\Vert$ be a norm on $\C^m$. Using the compactness of $X$ and holomorphic charts one finds  two constants $c>0$, $A>1$ and, for every point $\eta \in X$, a holomorphic map $\psi_\eta : B(\eta, 2c) \to \C^m$ such that
$$\frac{1}{A}d(\eta, \eta') \le \Vert \psi_\eta(\eta) - \psi_\eta(\eta')\Vert \le A d(\eta, \eta'),\;\; \forall \eta' \in B(\eta, 2c).
$$
For every $\eta \in X$, we set $J_\eta(\cdot):=d(\eta,\cdot)$, these functions $J_\eta$ clearly satisfy the assumptions $i)$ and $ii)$ of Lemma \ref{NewZalc}. Let us now consider a sequence of holomorphic maps  $f_n :\D \to X$ which is not normal at the origin $0$ of $\D$; we may assume that the $f_n$ are actually defined on the closure of $\D$. The classical Schwarz lemma, applied to the coordinates of $(\psi_\eta(\eta) - \psi_\eta(\eta'))$, shows that that the maps $f_n$ satisfy the Schwarz-type property $S(J,\alpha^{\pm},c, s)$ on $\D$ where $\alpha^+=\alpha^- =\alpha \in ]0,1[$ and $s(x):=A^2c\alpha x$. This corresponds to the assumption $1.$ of Lemma \ref{NewZalc}, while the assumption $2.$ immediately follows from the lack of equicontinuity,  on any neighbourhood of $0$,  of the sequence $(f_n)_n$. As we already observed, the
limits of  the renormalized sequence $(g_n)_n$ given by Lemma \ref{NewZalc} yield non constant entire  holomorphic curves in $X$.\\

 Miniowitz Lemma \cite{Min}  extends  Zalcman result to the setting of  quasi-regular  (\emph{qr}) maps  from the unit ball $\B^m$ of $\R^m$ to the Alexandroff compactification $\bar{\R}^m$ of $\R^m$.
 \begin{lemma}{\bf (Miniowitz)}
Let  $(f_n)_n$ be a sequence of  $K$-$qr$ maps from $\B^m$ to $\bar{\R}^m$ which is not normal at $0$. Then there exists a sequence $(x_n +r_n x)_n$
of affine contractions on $\R^m$, where  $r_n>0$ and $\lim r_n=0$,  such that after taking a subsequence, the renormalized sequence  
$(f_n(x_n+r_n x))_n$ is locally uniformly converging on $\R^m$ to a non constant $K$-$qr$ map  $\varphi : \R^m \to \bar{\R}^m$.
\end{lemma}
As we did  in the holomorphic case for Zalcman Lemma, one may deduce the above result  from our Lemma \ref{NewZalc}. To this end, one must replace $\D$ by $\R^m$ (see the remark at the end of the last section), and take $\bar{\R}^m$ endowed with the spherical distance as as metric space $X$. The key point is that $qr$ maps satisfy the following Schwarz Lemma (see \cite{Ric} Theorems 1.10).
 
 \begin{lemma}\label{LemSchwarzQR}
There exists a function  $s(m,K,\cdot)$ which is continuous and increasing from  $[0,1]$ to  $[0,1]$ such that $\vert f(x)\vert \le  s(m,K,\vert x\vert)$ for every $K$-quasi-regular  map $f:\B^m\to\B^m$
satisfying $f(0)=0$. 
\end{lemma}

\section{A Bloch principle for the Bedford-Pinchuk rescaling method in $\C^2$}\label{SecBP^2}

The rescaling method has been pioneered by Pinchuk (see \cite{Pin}). Roughly speaking, it consists of analysing the asymptotic behaviour of sequences of analytic objects $A_n$ converging to a boundary point of some domain $\Omega\subset \C^k$ from possible rigidity properties of  the limits of  $S_n(A_n)$, where the sequence $(S_n)_n$ is build to rescale the domain $\Omega$ to some simpler, but possibly unbounded, domain $D$. When $A_n$ is a sequence of holomorphic discs in $\Omega$, this method may yield to a sharp control  on the asymptotic behaviour of the
Kobayashi infinitesimal metric. When $A_n=f_n(\Omega)$ where the $f_n$ are automorphisms of 
$\Omega$, it may allow to characterize the domain $\Omega$ (see \cite{BP1}, \cite{BP2}, \cite{Be3}, \cite{V}, \cite{PSS}) . \\
Here we will consider the rescaling method in a setting corresponding to the seminal  works  \cite{BP1}, \cite{BP2} of
Bedford and Pinchuk: $\Omega$ will be  a domain in $\C^2$ whose boundary is smooth, pseudonconvex and of finite type
near $(0,0)$ and  $(A_n)_n$ will converge to $(0,0)$.\\

{\bf Geometry of finite type real hypersurfaces.} 
 A real hypersurface, defined in a neighbourhood of the origin in $\C^2$ by $\{ \rho=0\}$,  is said to be of type $m\in \N$ at $(0,0)$ if there exists an $m$-homogeneous real polynomial $H_m(z,\bar z)$ without harmonic terms, such that the defining function $\rho$ can be written in the following form in well chosen local coordinates:
\begin{eqnarray}\label{loccoor}
\rho(w_1,w_2)=2\Re\; w_2 + H_m(w_1,\bar{w}_1) + O\left(\vert w_1\vert^{m+1} +\vert w_1\vert \vert w_2\vert\right).
\end{eqnarray}
The hypersurface is pseudoconvex at $(0,0)$ if and only if the polynomial $H_m$ is subharmonic and, in that case, $m$ has to be even. We will suppose $\{\rho=0\}$
 globally pseudoconvex, of finite type $m$ at $(0,0)$, and will now describe some tools introduced by Catlin \cite{Cat}  to study such hypersurfaces.

We endow $\C^2$ with the norm $\Vert(w_1,w_2)\Vert_\infty:=\max(\vert w_1\vert, \vert w_2\vert)$, and the vector space of polynomials of degree at most $m$  in $w_1$, $\bar{w}_1$ with the analogous norm also denoted $\Vert\;\Vert_\infty$. We also fix a sufficently small  neighbourhood $U_0$ of $(0,0)$ in $\C^2$ and assume that the normal direction to $\{\rho =0\}$ at point $(0,0)$ is given by $(0,1)$. 
For every $\eta\in U_0$, we define the point $\hat \eta$ and the real (and positive when $\eta\in U_0^-$) quantity $\epsilon(\eta)$ by
$$\eta + (0,\epsilon(\eta)) =:\hat \eta \in \{\rho=0\}.$$ 

For every point $\eta \in U_0$ there exists a unique holomorphic automorphism $(\phi_\eta)^{-1}$ of $\C^2$ of the form 
$$\left(\phi_\eta\right)^{-1}(w_1,w_2):=\eta+(w_1, d_0(\eta)w_2 + \sum_{k=1}^m d_k(\eta) w_1^k)$$
such that
\begin{eqnarray}\label{rhocircphi}
\rho \circ \left(\phi_\eta\right)^{-1}-\rho(\eta)=
2\Re\; w_2 + 
\sum_{\stackrel{j+k\le m}{j,k>0}}
 a_{j,k}(\eta) w_1^j \bar{w}_1^k\ +O\left(\vert w_1\vert^{m+1} +\vert w_1\vert \vert w_2\vert\right).
\end{eqnarray}

The automorphisms $\phi_\eta$ depend continuously on $\eta$ and  we may assume that $\phi_{(0,0)}=\id$.
One thus easily sees   that $\phi_{\hat \eta}(\eta)=(0,-\tilde{\epsilon}(\eta))$ where $\tilde{\epsilon}(\eta)\sim \epsilon(\eta)\;\textrm{when}\; \eta\to (0,0).$ In particular, taking $U_0$ small enough we have
\begin{eqnarray}\label{espi/epsitil}
\phi_{\hat \eta}(\eta)=(0,-\tilde{\epsilon}(\eta)),\;\textrm{where}\;\vert \tilde{\epsilon}(\eta)\vert \le 2\vert\epsilon(\eta)\vert,\;\;\forall \eta\in U_0.
\end{eqnarray}
Next, for every $\eta\in U_0$ and every $\epsilon>0$ one defines the polynomial $P_\eta$ and the real quantity $\tau(\eta,\epsilon)$ by:
\begin{eqnarray}\label{tau}
P_\eta:=\sum_{\stackrel{j+k\le m}{j,k>0}} a_{j,k}(\eta) w_1^j \bar{w}_1^k,\;\; \;\;\tau(\eta,\epsilon):=\inf\{\tau >0\;\colon\; \Vert \frac{1}{\epsilon} P_\eta(\tau w_1,\tau\bar{w}_1)\Vert_\infty =1\}.
\end{eqnarray}
The functions $\tau(\eta,\cdot)$ are clearly increasing, and there actually exists $C_0, C_1>0$ such that the following estimates occur for $\eta\in U_0$ and $0<\epsilon\le\epsilon'\le 1$:
\begin{eqnarray}\label{eps/esp'}
C_0\epsilon^{\frac{1}{2}} \le \left(\frac{\epsilon}{\epsilon'}\right)^{\frac{1}{2}} \tau(\eta,\epsilon') \le \tau(\eta,\epsilon) \le
 \left(\frac{\epsilon}{\epsilon'}\right)^{\frac{1}{m}} \tau(\eta,\epsilon')\le C_1  \epsilon^{\frac{1}{m}}.
\end{eqnarray}
We now define the anisotropic dilations  $\Delta_\eta^\epsilon$ on $\C^2$ and the pseudo-balls $Q[\eta,\epsilon]$ by:
$$\Delta_\eta^\epsilon(w_1,w_2):=\left(\frac{w_1}{\tau(\eta,\epsilon)}, \frac{w_2}{\epsilon}\right),$$
$$Q[\eta,\epsilon]:=(\phi_\eta)^{-1} \circ (\Delta_\eta^\epsilon)^{-1} \left(\D^2\right)=(\phi_\eta)^{-1}\left(\D_{\tau(\eta,\epsilon)}\times \D_\epsilon\right).$$

The fundamental  properties of these pseudo-balls,  highlighted by Catlin (see \cite{Cat}), are the following.
\begin{proposition}\label{MainQ}
 If $U_0$ and $\alpha_0$ are small enough, then there exists constants $C_2,C_3,C_4 >1$ such that, for any $\eta,\eta'\in U_0$ and any $0<\epsilon<\alpha_0$, the following estimates occur if $\eta\in Q[\eta',\epsilon]$:
\begin{eqnarray*}
 (i)\;\epsilon(\eta)\le C_2(\epsilon(\eta') +\epsilon),\;\;
(ii)\;Q[\eta,\epsilon]\subset Q[\eta',C_3\epsilon],\;\;
(iii)\;Q[\eta',\epsilon]\subset Q[\eta,C_3\epsilon],\\
(iv)\; \frac{1}{C_4} \tau(\eta,\epsilon) \le \tau(\eta',\epsilon)\le C_4  \tau(\eta,\epsilon).
\end{eqnarray*}
\end{proposition}

We now introduce a family of functionals $(J_\eta)_{\eta \in U_0}$ on $\C^2$ by setting:
\begin{eqnarray}\label{Jeta}
J_\eta(w):=\Vert \Delta_{\hat{\eta}} ^{\epsilon(\eta)}\circ \phi_{\hat{\eta}}(w) - \Delta_{\hat{\eta}} ^{\epsilon(\eta)}\circ \phi_{\hat{\eta}}(\eta)\Vert_\infty.
\end{eqnarray}
The next proposition summarizes how  the finite type geometry will be used in the remaining of the paper, it crucially relies on the estimates given by Proposition \ref{MainQ}.
\begin{proposition}\label{LemQ}
If $U_0$ and $0<\epsilon_0<1$ are small enough, there exists a constant
$C_5>1$ such that:
\begin{itemize}
\item[i)] if $0<\epsilon<\epsilon_0$ and $\eta, \eta'\in U_0$ then: $ \eta' \in Q[{\hat \eta},\epsilon]\Rightarrow Q[{\hat \eta'},3\epsilon(\eta')]\subset
Q[{\hat \eta},C_5\epsilon]$.
\item[ii)] If $q_0,q_1,\cdots,q_p$ is a chain of points in $U_0^-$ such that
$3C_5^{p-1}\epsilon(q_0)<\epsilon_0$ and $J_{q_i}(q_{i+1}) <1$ for $0\le i\le p-1$,  then 
$q_p\in Q[\hat{q}_0, C(p)\epsilon(q_0)]$. In particular
$
\Delta_{{\hat q}_0}^{\epsilon(q_0)} \circ \phi_{{\hat q}_0} (q_p)\in C(p) \D^2,
$
where $C(p):={3C_5^{p-1}}$.
\end{itemize} 
\end{proposition}

\pro Let us prove  assertion i). By Proposition \ref{MainQ} i) we have $\epsilon(\eta')\le C_2(\epsilon(\hat \eta) + \epsilon)=C_2\epsilon$. Since $\hat\eta'=\eta' +(0,\epsilon(\eta'))$, we  get from the explicit
form of the automorphism $\phi_{\hat{\eta}}$ that $\phi_{\hat \eta}(\hat\eta')=\phi_{\hat \eta}(\eta')+(0,\frac{1}{d_0(\hat\eta)}\epsilon(\eta'))$. After shrinking $U_0$ we may assume that $\vert d_0(\hat \eta)\vert \ge \frac{1}{2}$ (recall that $\phi_{(0,0)}=Id$) and, since   $\eta' \in Q[{\hat \eta},\epsilon]$, we thus have
$$\phi_{\hat \eta}(\hat\eta') \in \D_{\tau(\hat\eta,\epsilon)} \times \D_\epsilon + (0,\frac{1}{d_0(\hat\eta)}\epsilon(\eta')) \subset  \D_{\tau(\hat\eta,\epsilon)} \times \D_{\epsilon+2\epsilon(\eta')} \subset  \D_{\tau(\hat\eta,K\epsilon)} \times \D_{K\epsilon},$$
where we have set $K:=1+2C_2$. This means that $\hat\eta'\in Q[{\hat \eta},K\epsilon]$.
Now, if $2K\epsilon_0 <\alpha_0$, Proposition \ref{MainQ} ii) implies that $Q[{\hat \eta'},2K\epsilon] \subset Q[{\hat \eta},2KC_3\epsilon]$. Finally, since
$3\epsilon(\eta') \le 3C_2\epsilon \le 2K\epsilon$, we get  $Q[{\hat \eta'},3\epsilon(\eta')] \subset  Q[{\hat \eta},C_5\epsilon]$ where $C_5:=2KC_3=(2+4C_2)C_3$.\\

Let us now prove the second assertion. We first observe that
\begin{eqnarray}\label{JQ}
J_\eta(w) <1 \Rightarrow w\in Q[\hat\eta, 3\epsilon(\eta)].
\end{eqnarray}
 Indeed, by definition, $J_\eta(w) <1$ implies that
\begin{eqnarray*}
\phi_{\hat \eta}(w)\in \phi_{\hat \eta}(\eta) + (\Delta_{\hat \eta}^{\epsilon(\eta)})^{-1}(\D^2)=
(0,-\tilde{\epsilon}(\eta))+\D_{\tau({\hat\eta},\epsilon(\eta))}\times \D_{\epsilon(\eta)}\subset \D_{\tau({\hat\eta},3\epsilon(\eta))}\times \D_{3\epsilon(\eta)},
\end{eqnarray*}
where  the last inclusion comes from (\ref{espi/epsitil}).
By  (\ref{JQ}) one has $q_i\in Q\left[\hat{q}_{i-1},3\epsilon(q_{i-1})\right]$ for every $1\le i\le p$ and, since
$3C_5^i \epsilon(q_0) \le \epsilon_0$  for $1\le i\le p-1$, we may use the first assertion inductively to  get
\begin{eqnarray*}
q_1&\in& Q\left[\hat{q}_0,3\epsilon(q_0)\right]\\
q_2&\in& Q\left[\hat{q}_1,3\epsilon(q_1)\right]\subset Q\left[\hat{q}_0,3C_5\epsilon(q_0)\right]\\
q_3&\in& Q\left[\hat{q}_2,3\epsilon(q_2)\right]\subset Q\left[\hat{q}_0,3C_5^2\epsilon(q_0)\right]\\
&\vdots&\\
q_{p-1}&\in& Q\left[\hat{q}_{p-2},3\epsilon(q_{p-2})\right]\subset Q\left[\hat{q}_0,3C_5^{p-2}\epsilon(q_0)\right]\\
q_{p}&\in& Q\left[\hat{q}_{p-1},3\epsilon(q_{p-1})\right]\subset Q\left[\hat{q}_0,3C_5^{p-1}\epsilon(q_0)\right].
\end{eqnarray*}
To conclude, one observe that  the definition of pseudo-balls, the estimates (\ref{eps/esp'}) and the fact that  $C(p):=3C_5^{p-1} \epsilon(q_0) \le \epsilon_0<1$ imply that
\begin{eqnarray*}
\Delta_{\hat q_0}^{\epsilon(q_0)} \circ \phi_{\hat q_0} \left(Q[\hat q_0,C(p)\epsilon(q_0)]\right)
\subset \frac{\tau(\hat q_0,C(p)\epsilon(q_0))}{\tau(\hat q_0,\epsilon(q_0))}\D\times C(p)\D \subset \sqrt{C(p)} \D\times C(p) \D \subset C(p)\D^2,
\end{eqnarray*}
 and thus get 
$ \Delta_{{\hat q}_0}^{\epsilon(q_0)} \circ \phi_{{\hat q}_0} (q_p)\in C(p) \D^2$.
\qed\\

{\bf The Bedford-Pinchuk rescaling method.} Let $\Omega$  be  a domain in $\C^2$ whose boundary is smooth, pseudonconvex and of finite type
near $(0,0)$ and  $(A_n)_n$ be a sequence of analytic objects which is converging  $(0,0)$ in $\Omega$. Let $U_0$ is a sufficently small ball centered at $(0,0)$  in $\C^2$, then we may assume that  $\Omega \cap U_0$ is a domain of the form $U_0^-=\{ w\in U_0\;\colon\; \rho(w)<0\}$ like those that we just studied. Note that  $U_0^-$ is pseudoconvex.\\

Let $(\eta_n)_n$ be a sequence converging to $(0,0)$ in $U_0^-$. When $n$ is big enough, the point $\eta_n$
is sufficently close to $(0,0)$ and thus there exists $\epsilon_n:=\epsilon(\eta_n)>0$ such that
\begin{eqnarray}\label{eps(eta)}
\eta_n + (0,\epsilon_n) =:\hat{\eta}_n \in \{\rho=0\}. 
\end{eqnarray}
One sees from (\ref{rhocircphi}), (\ref{tau}), and (\ref{eps/esp'}) that
\begin{eqnarray}\label{rhodilat}
\rho \circ  \left(\phi_{\hat{\eta}_n}\right)^{-1} \circ \left(\Delta_{\hat{\eta}_n}^{\epsilon_n}\right)^{-1}=2\Re\; w_2 + P_n(w_1,\bar{w}_1) + O(\tau({\hat \eta}_n, \epsilon_n)),
\end{eqnarray}
where $P_n$ is the polynomial of degree at most $m$ satisfying $\Vert P_n \Vert_\infty =1$ given by 
$$P_n(w_1,\bar{w}_1):=\frac{1}{\epsilon_n} P_{\hat{\eta}_n}(\tau({\hat \eta}_n, \epsilon_n) w_1,\tau({\hat \eta}_n, \epsilon_n)\bar{w}_1).$$
Consider the sequence $(S_n)_n$ of rescaling automorphisms of $\C^2$   defined by $S_n := \Delta_{\hat{\eta}_n}^{\epsilon_n}
\circ \phi_{\hat{\eta}_n}.$ Then
\begin{eqnarray}\label{Pinc1} 
S_n : U_0^- \to U_n^-:=S_n(U_0^-),\;\;S_n(\hat{\eta}_n)=(0,0)\;\;\textrm{and}\;\; \lim_n S_n({\eta}_n)=(0,-1).
\end{eqnarray}
After taking a subsequence, we may assume that the sequence of polynomials $(P_n)_n$ is converging to some polynomial $P$ such that $\Vert P\Vert_\infty =1$. Then, it follows from (\ref{rhodilat}) that the sequence of bounded domains $(U_n^-)_n$ is converging in the Hausdorff sense to some unbounded rigid polynomial domain $D_P$ associated to $P$:
\begin{eqnarray}\label{Pinc2}
U_n^- \to \{(w_1,w_2)\in \C^2\;\colon\; 2\Re\;w_2 + P(w_1,\bar{w}_1) <0\}=:D_P.
\end{eqnarray}
The domains $U_n$ are pseudoconvex and, using the tomato-can principle, one sees that the limit
domain $D_P$ is pseudoconvex too. The polynomial $P$ is therefore subharmonic and, like the 
$P_n$, does not contain any constant or  harmonic term.\\

The rescaling method requires to produce  limits of $(S_n(A_n))_n$. Although the limit domain $D_P$ is taut, this is not obvious
because $S_n(A_n)$ is generally not contained in $D_P$ or in any similar domain. However, our next theorem explains why this is always  possible.\\

{\bf  A general normality statement.}

\begin{theorem}\label{TheoNormaPin}
Let $\Omega$ be a domain in $\C^2$ whose boundary is smooth, pseudonconvex and of finite type $m$
near $(0,0)$. Let $\omega$ be a domain in $\C^k$ and $(f_n)_n$ be a sequence of holomorphic maps
from $\omega$ to $\Omega$ such that $\lim_n f_n(a_0)=(0,0)$ for some fixed point $a_0\in \omega$.

Let $(S_n)_n$ be the sequence of automorphisms of $\C^2$ associated by the above described rescaling method
to the sequence of points $(\eta_n)_n:=(f_n(a_0))_n$.

Then $(S_n\circ f_n)_n$ is normal and its limits are  holomorphic maps from
$\omega$ to some rigid  domain $D_P:=\{(w_1,w_2)\in \C^2\;\colon\; 2\Re\; w_2 + P(w_1,{\bar w}_1) <0\}$ where $P$ is a subharmonic and non harmonic polynomial of degree at most $m$.
\end{theorem}
As we shall see, this result  once again illustrates the {\it Bloch principle : as all entire holomorphic curve in the limit domain $D_P=\lim_n S_n(\Omega)$ are constant,  the sequence $(S_n \circ f_n)_n$ is normal for any sequence $(f_n)_n$ of holomorphic discs in $\Omega$ whose centers $f_n(0)$ converge to $(0,0)$.}
 
The proof we gave in \cite{Be} was based on some renormalization technique very similar to Zalcman's one (see Lemma 3.1 in \cite{Be} and required some
delicate integration arguments for specific pseudo-metrics.

The proof we present here avoids all these technical difficulties, it entirely relies on our \ref{NewZalc} lemma combined with a version of Bloch's principle. The Lemma \ref{NewZalc} will be applied 
 with $X=U_0^-$ endowed by the metric $d$ induced by $\Vert\;\Vert_\infty$, the functionals $J_\eta$ being defined by (\ref{Jeta}).
The required assumptions are clearly satisfied, in particular the classical Schwarz lemma implies that any $f\in{\mathcal O}(\D, U_0^-)$ satisfies the
property $S(J,1,1,s)$ with $s(u)=\vert u\vert$.  Note that $\alpha^{\pm}=c=1$. 
From now on, for any map $\varphi$ to $U_0^-$,  we shall denote by $\hat{\varphi}(t)$ the point $\widehat{\varphi(t)}$.\\

The heart of the proof of Theorem  \ref{TheoNormaPin} lies in Lemma \ref{LemPBA} below.
As we will explain in the next section, this lemma implicitely contains
Catlin's estimates on the Kobayashi infinitesimal metric.\\

\begin{lemma}\label{LemPBA} For every $0<k\le 1$ there exists $0<r_0<1$ and $c>0$ such that
$$\left( f\in {\mathcal O}(\D, U_0^-)\;\;\textrm{and }\;\; \Vert f(0)\Vert_\infty \le c\right)\Rightarrow J_{f(0)}(f(t)) < k,\;\forall t\in \D_{r_0}.$$
\end{lemma}
\pro
We proceed by contradiction. If this is not true, we may find $k>0$, a sequence $(\epsilon'_n)_n$ converging to $0$ in 
$\D$ and a sequence  $(f_n)_n$  in ${\mathcal O}(\D, U_0^-)$ such that $\lim_n f_n(0)=(0,0)$ and  $J_{f_n(0)}(f_n(\epsilon'_n)) \ge k$ for every $n\in \N$.

As the boundary of $U_0^-$ does not contain analytic discs through $(0,0)$, it follows from Montel's theorem that $(f_n)_n$ is locally uniformly converging to $(0,0)$ and thus,
after replacing $f_n(z)$ by $f_n(\frac{z}{2})$ and $\epsilon'_n$ by $2\epsilon'_n$, we may assume that $f_n\in {\mathcal O}(\D_2,U_0^-)$. 
Then Lemma \ref{NewZalc}, applied with $t'_n=0$, provides a sequence of affine contractions $(r_n)_n$ on $\C$, with $r_n(0)\in \D$,  and a sequence of positive numbers $(R_n)_n$ converging to $+\infty$ such that  the maps $g_n:=f_n\circ r_n$ are defined on $\D_{R_{n}+1}$ and 
\begin{eqnarray}\label{tt'Q}
J_{g_n(0)}\left(g_n(1)\right) \ge k\;\;\textrm{and}\; J_{g_n(t)}\left(g_n(t+u)\right)\le \vert u\vert,\; \forall t\in \D_{R_n}\;\forall u\in \D.
\end{eqnarray}

As $g_n(0)=f_n(r_n(0))$  and $\vert r_n(0)\vert <1$,  one has $\lim_n g_n(0)=(0,0)$. Taking subsequences, we may also assume that $R_n> n$ and $3 C_5^n \epsilon(g_n(0)) \le \epsilon_0$ for every $n\in \N$, where $C_5$ and $\epsilon_0$ are the constants given by Proposition \ref{LemQ}. 

For  $p\le n\in\N$ and $t\in \D_p$, let $0=:t_0, t_1, t_2,\cdots, t_{p-1}, t_p:=t$ be a chain of points in $\D_p$ such that
$\vert t_i-t_{i-1}\vert <1$ for $1\le i\le p$. 
According to (\ref{tt'Q}) we may apply the second assertion of Proposition \ref{LemQ} to the chain of points $\left(q_i:=g_n(t_i)\right)_{0\le i\le p}$ and get
\begin{eqnarray}\label{majMont}
G_n(t):=\Delta_{{\hat g}_n(0)}^{\epsilon(g_n(0))} \circ \phi_{{\hat g}_n(0)} \circ g_n(t) \in C(p) \D^2,\;\;\textrm{for every}\;\; t\in D_p.
\end{eqnarray}

 By (\ref{majMont}) and Montel's theorem, the sequence $(G_n)_n$ is locally uniformly converging on $\C$,
after taking a subsequence. Setting $\eta_n:=g_n(0)$ one sees that  $G_n=S_n\circ g_n$ and then  (\ref{Pinc1}) and (\ref{Pinc2})
show that the  limit $G$ is a holomorphic entire curve in $\overline{D_P}$ and $G(0)=(0,-1)$. The maximum modulus principle applied to the negative subharmonic function $\rho(t):= \Re\; G_2(t) +P(G_1(t))$, which satisfies $\rho(0)=-1$, shows that the curve $G$ is actually contained in $D_P$.
Finally, since $$\Vert G_n(0)-G_n(1)\Vert_\infty= J_{g_n(0)}(g_n(1)) \ge k>0,$$  $G$ is not constant. This is the expected contradiction (see Lemma \ref{BroRig} below)
.\qed\\

To prove Theorem \ref{TheoNormaPin} we first need is a localization statement for the sequence $(f_n)_n$.
\begin{lemma}\label{LemLoc}
The sequence $(f_n)$ is locally uniformly converging to $(0,0)$ on $\omega$.
\end{lemma}

\pro When $\Omega$ bounded, this follows immediately from Montel's theorem and the fact that there are no analytic disc in $b\Omega$ through $(0,0)$. When $\Omega$ is not bounded, this is still true but the proof is more involved (see \cite{Be2} Proposition 2.1 and \cite{FS} for the existence of $p.s.h$ peak functions, or \cite{Bh} Result 3.2, \cite{TT} Proposition 2.2).
\qed\\

We may now end the {\it proof of Theorem \ref{TheoNormaPin}}. We endow $\C^k$ with the euclidean norm $\vert\cdot\vert$. Let $a\in \omega$, we have to show that
the sequence $(S_n\circ f_n)_n$ is normal on some neighbourhood of $a$. Let us recall that $S_n=\Delta_{{\hat f}_n(a_0)}^{\epsilon(f_n(a_0))} \circ \phi_{{\hat f}_n(a_0)}$. 

Let $\gamma:[0,1] \to \omega$ be a continuous path connecting $a_0$ to $a$ in $\omega$ and let $r:=d(\gamma([0,1]), b\omega)$.
Let $C_5>0, \epsilon_0, r_0, c$ be the constants  given by  Proposition \ref{LemQ},  and Lemma \ref{LemPBA} applied with $k=1$. 
Let $a'\in B(a,r_0r)\setminus\{a\}$. We take $0=t_0<t_1<t_2<\cdots<t_{p-1}=1$ such that $0<\vert \gamma(t_i)-\gamma(t_{i-1})\vert < r_0 r$ for $1\le i\le p-1$, and set $a_i:=\gamma(t_i)$ for $0\le i\le p-1$ and $a_p:=a'$.
We both have $B(a_i,r)\subset \omega$ for $0\le i\le p-1$ and $\vert a_i-a_{i-1}\vert <r_0r$ for $1\le i\le p$.

Taking $n\ge n_0$ with $n_0$ big enough will guarantee that $3C_5^{p} \epsilon(f_n(a_0)) < \epsilon_0$ and, by Lemma \ref{LemLoc}, that $f_n\left(\cup_{0\le i\le p-1} B(a_i,r)\right)\subset U_0^-$ and $\max_{0\le i\le p} \Vert f_n(a_i))\Vert_\infty < c$. Then, applying Lemma \ref{LemPBA} to the maps
$\D\ni u \mapsto f_n\left( a_{i-1}+\frac{ru}{\vert a_i-a_{i-1}\vert}(a_i-a_{i-1})\right)$ yields $J_{f_n(a_{i-1})}\left(f_n(a_i)\right) < 1$ for $1\le i\le p$ and therefore, 
the second assertion of Proposition \ref{LemQ}  applied to the chain of points $q_i=f_n(a_i)$ yields 
\begin{eqnarray*}
S_n\circ f_n(a')=\Delta_{{\hat f}_n(a_0)}^{\epsilon(f_n(a_0))} \circ \phi_{{\hat f}_n(a_0)} \circ f_n(a') \in C(p) \D^2,\;\;\forall a'\in B(a,r_0r).
\end{eqnarray*}
The conclusion follows from Montel's theorem.\qed\\
\begin{lemma}\label{BroRig}
If $P$ is a subharmonic and non harmonic function on $\C$ then the domain $D_P:=\{(w_1,w_2)\in \C^2\;\colon\; \Re\; w_2 + P(w_1) <0\}$ is Brody-hyperbolic.
\end{lemma}
\pro If $f:\C\to D_P$ is holomorphic then $\Re\; f_2 +P\circ f_1$ is subharmonic negative, and therefore constant, on $\C$. The function $P\circ f_1$
is thus harmonic on $\C$ and therefore $f_1$ is constant on $\C$. It follows that $\Re f_2$ and $f_2$ are constant too.\qed\\

\section{Estimates for the Kobayashi pseudo-metric near finite type boundaries in $\C^2$}

Sharp estimates of the Kobayashi pseudo-metric near finite-type boundaries in $\C^2$ were first obtained  by Catlin \cite{Cat}.  His  proof, consisting of minimizing the Carathéodory metric, gives a stronger result but is technically difficult (bumping, $L^2$-estimates of H\"{o}rmander). An alternative approach, based on a delicate construction of $p.s.h$-barriers due to Fornaess and Sibony, was carried out by Coupet and Sukhov in \cite{CS}.

The proof we present here, which  is essentially based on the Lemma \ref{LemPBA},  is absolutely elementary.\\

Like in the former section we consider   a domain  $\Omega$  in $\C^2$, possibly unbounded,  whose boundary is smooth, pseudonconvex and of finite type
near $(0,0)$.  If  $U_0$ is a sufficently small ball centered at $(0,0)$,  then $U_0^-:=\Omega \cap U_0$ is a bounded pseudo-convex domain and there exists a smooth function $\rho$ such that  $U_0^-=\{ w\in U_0\;\colon\; \rho(w)<0\}$.

As we saw at the beginning of section \ref{SecBP^2}, to every point $\eta \in U_0^-$ are attached an automorphism  $\phi_\eta$
of $\C^2$ and positive quantities $\epsilon(\eta)$ and $\tau(\eta,\epsilon(\eta))$. Let us recall that $\epsilon(\eta)$ is the distance
between $\eta$ and the piece of hypersurface $\{\rho = 0\}$ in the direction of its gradient at $(0,0)$, and that we may assume that
$\phi_\eta =id$ when $\eta=(0,0)$. This last condition being achieved after replacing $\Omega$ by $\phi_{(0,0)}(\Omega)$.

Using these quantities, we may define a pseudo-metric $M_{U_0^-}$ on $U_0^-$ by setting
$$M_{U_0^-} (\eta,X):=\max\left(\frac{\vert (\phi'_\eta(\eta)\cdot X)_1\vert}{\tau(\eta,\epsilon(\eta))},\frac{\vert (\phi'_\eta(\eta)\cdot X)_2\vert}{\epsilon(\eta)}\right)$$
for every $\eta \in U_0^-$ and every $X\in \C^2$.
Our aim is to show that the Kobayashi pseudo-metric $K_\Omega$ of the domain $\Omega$ is equivalent to $M_{U_0^-}$
on $U_0^-$, that is:
$$K_\Omega(\eta, X) \approx  M_{U_0^-}(\eta,X),\;\;\forall \eta \in U_0^-,\; \forall X\in \C^2.$$

Let us recall that the Kobayashi pseudo-metric is defined by:
$$K_\Omega(\eta,X):=\inf \{\frac{1}{R}\;\colon\; \exists f\in {\mathcal O}(\D,\Omega)\;\textrm{such that}\; f(0)=\eta\;\textrm{and}\; f'(0)=R X\}.$$

As the injection $i :U_0^- \to \Omega$ decreases the Kobayashi metrics we have $K_\Omega \le K_{U_0^-}$, and  it easily follows from Lemma \ref{LemLoc} that $K_{U_0^-} \le 2 K_\Omega$ after maybe reducing $U_0^-$.
Thus

\begin{eqnarray}\label{LocK}
K_\Omega(\eta, X)\le K_{U_0^-}(\eta,X) \le 2 K_\Omega(\eta,X),\;\;\forall \eta \in U_0^-,\; \forall X\in \C^2.
\end{eqnarray}
On the other hand, it follows from the very definitions that
$$ K_{U_0^-}(\eta,X) \lesssim M_{U_0^-}(\eta,X),\;\;\forall \eta \in U_0^-,\; \forall X\in \C^2.$$
It is therefore sufficent  to show  that  $ M_{U_0^-}(\eta,X)\lesssim  K_{U_0^-}(\eta,X)$ for $\eta \in U_0^-$ and $X\in \C^2$.
As we shall see, this is essentially what says Lemma \ref{LemPBA}.\\

Let $f\in {\mathcal O}(\D, U_0^-)$ such that $f(0)=\eta$ and $f'(0)=RX$ where $R>0$.
By Lemma \ref{LemPBA}, we have $J_{\eta}(f(t)) \le 1$ for every $t \in \D_{r_0}$, which, by definition of $J$, means that
$$\Vert \Delta_{\hat \eta}^{\epsilon(\eta)} \circ \phi_{\hat \eta}(f(t)) -  \Delta_{\hat \eta}^{\epsilon(\eta)} \circ \phi_{\hat \eta}(f(0))\Vert_\infty \le 1,\;\;\forall t\in \D_{r_0}.$$
By the Schwarz lemma, this implies that
\begin{eqnarray}\label{EstimKobaBrut}
\frac{1}{R}\ge  r_0 \Vert  \Delta_{\hat \eta}^{\epsilon(\eta)} \circ \phi'_{\hat \eta}(\eta) \cdot X \Vert_\infty.
\end{eqnarray}
To conclude it remains to compare $M_{U_0^-}(\eta,X)=
\Vert  \Delta_{ \eta}^{\epsilon(\eta)} \circ \phi'_{ \eta}(\eta) \cdot X \Vert_\infty$ with  the right hand side of (\ref{EstimKobaBrut})  
and therefore to check that the two linear automorphisms  $ \Delta_{\hat \eta}^{\epsilon(\eta)}\circ  (\Delta_{ \eta}^{\epsilon(\eta)})^{-1}$ and $\phi'_{\hat \eta}(\eta)\circ (\phi'_{ \eta}(\eta))^{-1}$,
as well as their inverses, are bounded in norm. 

For $\phi'_{\hat \eta}(\eta)\circ (\phi'_{ \eta}(\eta))^{-1}$, it suffices  to reduce $U_0$ so that $\phi_\eta$ and $\phi_{\hat\eta}$ are sufficently close to $\phi_{(0,0)}=id$. For
$ \Delta_{\hat \eta}^{\epsilon(\eta)}\circ  (\Delta_{ \eta}^{\epsilon(\eta)})^{-1}=\frac{\tau(\eta,\epsilon(\eta))}{\tau({\hat\eta},\epsilon(\eta))} id \times id$, we first observe that (\ref{espi/epsitil}) tells us that
$\phi_{\hat \eta}(\eta) =(0,-\tilde\epsilon(\eta)) \in \D_{\tau(\hat \eta,2\epsilon(\eta))}\times \D_{2\epsilon(\eta)}$,
which means that $\eta\in Q[\hat\eta,2\epsilon(\eta)]$. Then, by Proposition \ref{MainQ} (iv), we have
$\frac{1}{C_4} \le \frac{\tau( \eta,2\epsilon(\eta))}{\tau(\hat \eta,2\epsilon(\eta))}\le C_4$ which finally, using
(\ref{eps/esp'}), yields $( \frac{1}{2})^{\frac{1}{2}}\frac{1}{C_4}\le \frac{\tau(\eta,\epsilon(\eta))}{\tau({\hat\eta},\epsilon(\eta))} \le (\frac{1}{2})^{\frac{1}{m}}C_4$.
This shows that, for $U_0$ small enough, (\ref{EstimKobaBrut}) implies
\begin{eqnarray}\label{EstimKobaBrut2}
M_{U_0^-}(\eta,X)\le    A(r_0,C_4)  K_{U_0^-}(\eta,X),\;\;\forall \eta \in U_0^-,\; \forall X\in \C^2
\end{eqnarray}
 where  $A(r_0,C_4)$ is a function which only depends on the quantities $r_0$ (given by Lemma \ref{LemPBA}) and
 $C_4$ (given by Proposition \ref{MainQ}).
 Writing $M_{U_0^-}(\eta,X)$ more explicitely and taking (\ref{LocK}) into account,  (\ref{EstimKobaBrut2}) becomes
 $$\frac{1}{2A(r_0,C_4)} \max \left( \frac{\vert X_1\vert}{\tau(\eta,\epsilon(\eta))},\frac{\vert\rho_{w_1}(\eta) X_1 + \rho_{w_2}(\eta) X_2\vert}{\epsilon(\eta)}\right) \le K_{\Omega}(\eta,X),\;\;\forall \eta \in \Omega\cap U_0,\;\forall X\in \C^2.$$
 Finally, as $\epsilon(\eta) \sim \vert \rho(\eta)\vert$  and $\tau(\eta,\epsilon(\eta))\approx \tau(\eta,\vert \rho(\eta)\vert)$ when $\eta\to (0,0)$, we may express this estimates in terms of the defining function $\rho$ as follows.
 
\begin{theorem}\label{TheoKoba}
Let $\Omega$ be a domain in $\C^2$ whose boundary is smooth, pseudonconvex and of finite type $m$
near $(0,0)$. Let $U_0$ be s sufficently small ball centered at the origin in $\C^2$ and $\rho$ be a smooth defining  function for $b\Omega\cap U_0$.
Assume that $\rho_{w_1}^{(r)}(0,0)=0$ for $1\le r\le m$. Then there exists a constant $A>1$ such that the Kobayashi pseudo-metric $K_\Omega$ of $\Omega$ satisfies the 
following estimates:
$\frac{1}{A} \max \left( \frac{\vert X_1\vert}{\tau(\eta)},\frac{\vert\rho_{w_1}(\eta) X_1 + \rho_{w_2}(\eta) X_2\vert}{\vert\rho(\eta)\vert}\right) \le K_{\Omega}(\eta,X),\;\;\forall \eta \in \Omega\cap U_0,\;\forall X\in \C^2$
where   
$\tau(\eta):=\min \{ (\frac{(j+k)!\vert\rho(\eta)\vert}{\vert (\rho_{w_1}^{(j)})_{\bar w_1}^{(k)}(\eta)\vert })^{\frac{1}{j+k}}\; \colon\;j,k\ge 1,\; j+k\le m\}$.
\end{theorem}

\end{document}